\documentclass[a4paper,10pt]{article}
\pdfoutput=1
\usepackage{amsfonts, amsthm, amsmath, amssymb}
\usepackage{fancyhdr}
\usepackage{graphicx}

\DeclareMathOperator*{\sign}{sgn}
\DeclareMathOperator*{\BMO}{BMO}
\providecommand{\abs}[1]{\ensuremath{\left\lvert #1 \right\rvert}}
\providecommand{\norm}[1]{\ensuremath{\left\Vert #1 \right\Vert}}
\providecommand{\vv}[1]{\textquotedblleft #1\textquotedblright}
\providecommand{\floor}[1]{\ensuremath{\left[ #1 \right]}}
\newcommand{\confrac}[2]{
\frac{\displaystyle{
\strut\hfill{#1}\hfill\;\vrule}}
{\displaystyle{
 \strut\vrule\;\hfill{#2}\hfill}}}

\newtheorem{theo}{Theorem}
\newtheorem{lem}{Lemma}
\newtheorem{prop}{Proposition}

\theoremstyle{definition}

\newtheorem{rem}{Remark}

\newcommand{\Q}{\mathbb{Q}}
\newcommand{\R}{\mathbb{R}}

\newcommand{\Z}{\mathbb{Z}}

\title{Generalized Brjuno functions associated to $\alpha$-continued fractions}

\author{\small{\textsc{Laura Luzzi\footnote{Scuola Normale Superiore, Piazza dei Cavalieri 7, 56123 Pisa, Italy, e-mail: \texttt{l.luzzi@sns.it}}, Stefano Marmi\footnote{Scuola Normale Superiore, Piazza dei Cavalieri 7, 56123 Pisa, Italy, e-mail: \texttt{s.marmi@sns.it}},  Hitoshi Nakada\footnote{Department of Mathematics, Keio University, 3-14-1 Hiyoshi, Kohoku-ku, Yokohama, 223-8522, Japan, e-mail: \texttt{nakada@math.keio.ac.jp}}, and Rie Natsui\footnote{Department of Mathematics, Japan Women's University, 2-8-1 Mejirodai, Bunkyou-ku, Tokyo, 112-8681, Japan, e-mail: \texttt{natsui@fc.jwu.ac.jp}}}}}
\date{}
\begin{document}
\maketitle

\begin{abstract}
For $0 \leq \alpha \leq 1$ given, we consider the $\alpha$-continued fraction expansion of a real number obtained by iterating the map $$A_{\alpha}(x)=\abs{x^{-1}-\floor{x^{-1}+1-\alpha}}$$ defined on the interval $\mathbb{I}_{\alpha}=(0,\bar{\alpha})$, with $\bar{\alpha}=\max(\alpha,1-\alpha)$. These maps generalize the classical (Gauss) continued fraction map which corresponds to the choice $\alpha=1$, and include the nearest integer ($\alpha=1/2$) and by-excess ($\alpha=0$) continued fraction expansion. 
To each of these expansions and to each choice of a positive function $u$ on the interval $\mathbb{I}_\alpha$ we associate a generalized Brjuno function $B_{(\alpha,u)}(x)=\sum_{n=0}^\infty \beta_{n-1} u(x_{\alpha,n})$, where $x_{\alpha,n}=A_{\alpha}(x_{\alpha,n-1})$ for all $n \geq 1$, $x_{\alpha,0}=\abs{x-\floor{x+1-\alpha}}$, $\beta_{\alpha,n}=x_{\alpha,0} \cdots x_{\alpha,n}$, $\beta_{\alpha,-1}=1$. When $\alpha=1/2$ or $\alpha=1$, and $u(x)=-\log(x)$, these functions were introduced by Yoccoz in his work on linearization of holomorphic maps. Their regularity properties, including BMO regularity and their extension to the complex plane, have been thoroughly investigated. We compare the functions obtained with different values of $\alpha$ and we prove that the set of $(\alpha,u)$-Brjuno numbers does not depend on the choice of $\alpha$ provided that $\alpha \neq 0$.
We then consider the case $\alpha=0$, $u(x)=-\log(x)$ and we prove that $x$ is a Brjuno number (for $\alpha \neq 0$) if and only if both $x$ and $-x$ are Brjuno numbers for $\alpha=0$.   
\bigskip \par
\noindent\textbf{Subject Classification:} Primary: 11J70; Secondary: 37F50\\
\textbf{Keywords:} continued fractions, Brjuno function, approximations of real numbers
\end{abstract}

\thanks{}

\section*{Introduction}
Let $x \in \R \setminus \Q$ and let $\{\frac{p_n}{q_n}\}_{n\geq 0}$ be the sequence of the convergents of its continued fraction expansion. A \emph{Brjuno number} is an irrational number $x$ such that $\sum_{n=0}^{\infty} \frac{\log q_{n+1}}{q_n} < \infty$. 
Almost all real numbers are Brjuno numbers, since for all diophantine numbers one has $q_{n+1}=\mathbb{O}(q_n^{\tau+1})$ for some $\tau\geq 0$. But some Liouville numbers also verify the Brjuno condition, e.g. $l=\sum_{n=0}^{\infty} 10^{-n!}$. 
The importance of Brjuno numbers comes from the study of $1$-dimensional analytic small divisors problems. The most important results are due to Yoccoz \cite{Yo}, who proved that the Brjuno condition is optimal for the problem of linearization of germs of analytic diffeomorphisms with a fixed point (and also for linearizing analytic diffeomorphisms of the circle provided that they are sufficiently close to a rotation \cite{Yo2}.) 

The set of Brjuno numbers is invariant under the action of the modular group $PGL(2,\Z)$ and can be characterized as the set where the Brjuno function $B: \R \setminus \Q \to \R \cup \{\infty\}$ is finite, where 
$$ B(x)=-\sum_{n=0}^\infty \beta_{n-1} \log x_n,$$
and $\beta_{-1}=1$, $\beta_n=x_0 \cdots x_n$, $x_{n+1}=\{x^{-1}\}$. 
The Brjuno function is $\Z$-periodic and satisfies the functional equation
$$B(x)=-\log x + x B(x^{-1}), \quad x \in (0,1),$$
which allows to interpret $B$ as a \emph{cocycle} under the action of the modular group (see \cite{MMY1}, \cite{MMY2} for details) and to study its regularity properties. For example one can prove that $B \in \bigcap_{p \geq 1} L^p(0,1)$ and even that $B \in \BMO(\mathbb{T}^1)$ \cite{MMY}. 

The aim of this paper is to extend some of the results known for $B$ to more general objects, obtained replacing the logarithm with an arbitrary positive $C^{1}$ function on $(0,1)$, with a singularity at the origin, and using $\alpha$-continued fractions instead of the Gauss map.  
\bigskip \par
Let $0 \leq \alpha \leq 1$, $\bar{\alpha}=\max(\alpha,1-\alpha)$. The \emph{$\alpha$-continued fraction expansion} of a real number $x \in (0, \bar{\alpha})$ is associated to the iteration of the map $A_{\alpha} : (0,\bar{\alpha}) \to (0,\bar{\alpha})$,
$$A_{\alpha}(x)=\abs{\frac{1}{x}-\floor{\frac{1}{x}+1-\alpha}}$$
These maps were introduced by one of the authors \cite{Nak} and include the standard continued fraction map ($\alpha=1$), the nearest integer ($\alpha=1/2$) and by-excess continued fraction map ($\alpha=0$) as special cases. 
For all $\alpha \in (0,1]$ these maps are expanding and admit a unique absolutely continuous invariant probability measure $d\mu_{\alpha}=\rho_{\alpha}(x) dx$ \cite{Nak} \cite{MCM} \cite{LM} whose density is bounded from above and from below by a constant dependent on $\alpha$. In the case $\alpha=0$ there is an indifferent fixed point and $A_{\alpha}$ does not have a finite invariant density but it preserves the infinite measure $d\mu_0(x)=\frac{dx}{1-x}$. 

In \S 1, given a positive $C^1$ function $u$ on $(0,1)$ with a singularity at the origin, we define a generalized Brjuno function $B_{\alpha,u}(x)=\sum_{n=0}^\infty \beta_{\alpha,n-1} u(x_{\alpha,n})$, where $\beta_{\alpha,n}=x_{\alpha,0} \cdots x_{\alpha,n}$ for $n \geq 0$, $\beta_{\alpha,-1}=1$ and the sequence $x_{\alpha,n}$ is obtained by iterating the $\alpha$-continued fraction map: $x_{\alpha,n+1}=A_{\alpha}(x_{\alpha,n})$, $n \geq 0$. 
We then prove that the set of convergence of $B_{\alpha,u}$ does not depend on the choice of $\alpha$ provided that $\alpha>0$, and that actually the difference $B_{\alpha,u}-B_{\alpha',u}$ between any two functions is $L^{\infty}$. 
The generalized Brjuno functions also satisfy a functional equation under the action of the modular group. 

In \S 2 we investigate the relation between the standard Brjuno function $B_1$ (corresponding to the choice $\alpha=1, u(x)=-\log x$) and the Brjuno function $B_0$ obtained replacing the Gauss map with the by-excess continued fraction map (i.e. $\alpha=0$, $u(x)=-\log x$). We prove that an irrational number $x$ is a Brjuno number if and only if $B_0(x)< \infty$ and $B_0(-x)<\infty$. Moreover the difference between $B_1$ and the even part of $B_0$ is bounded and numerical simulations suggest that it is $\frac{1}{2}\,$-H\"older-continuous.   

\section{The $(\alpha,u)$-Brjuno functions}  
\subsection{Some basic notions}
Fix $\alpha$, $0 \leq \alpha \leq 1$, let $\bar{\alpha}=\max(\alpha,1-\alpha)$, and define the map $A_{\alpha}$ from 
${\mathbb I}_{\alpha} = (0, \bar{\alpha})$ onto itself as follows:
\[
   A_{\alpha}(x) = \abs{\frac{1}{x}- 
   \left[ \frac{1}{x} \right]_{\alpha}}
\]
for $x \in {\mathbb I}_{\alpha}$, where 
$\floor{x}_{\alpha} = \floor{x+1-\alpha}$, and $\floor{\;\,}$ denotes the integer part. 
We put
\[
a_{\alpha}(x) = \floor{\frac{1}{x}}_{\alpha}, 
\qquad \varepsilon_0(x) = \sign(x-\floor{x}_\alpha)
\]
and define 
\[
\left\{ 
\begin{array}{rcl} 
a_{\alpha, n}(x) & = & a_{\alpha}(A_{\alpha}^{n-1} (x)) \\
\varepsilon_{\alpha, n} & = & \sign\left( \frac{1}{x_{n-1}}- 
   \left[ \frac{1}{x_{n-1}} \right]_{\alpha}\right) .
\end{array} \right. 
\]
These maps are factors of the $\alpha$-continued fraction transformations introduced by H. Nakada \cite{Nak}. 
In particular $A_1$ and $A_0$ are respectively the standard continued fraction transformation and the continued fraction 
transformation associated to the \vv{$-$} or \vv{by excess} expansion.  
We also define the following matrices : 
\[
M_{\alpha , n} \, = \, 
\begin{pmatrix}
0  &  \varepsilon_{\alpha , n-1} \\
1  &  a_{\alpha , n} 
\end{pmatrix}
\]
for $n \ge 1$, 
\begin{equation} \label{convergents}
\begin{pmatrix}
p_{\alpha , n-1}  & p_{\alpha , n } \\ 
q_{\alpha , n-1}  & q_{\alpha , n}  \\ 
\end{pmatrix}
 =  
M_{\alpha , 1} \cdot M_{\alpha , 2} \cdots M_{\alpha , n} 
\end{equation}
for $n \ge 1$, and 
\begin{equation} \label{p0}
M_{\alpha , 0} \, = \, 
\begin{pmatrix}
p_{\alpha , -1}  & p_{\alpha , 0}  \\ 
q_{\alpha , -1}  & q_{\alpha , 0}  \\ 
\end{pmatrix}
 =  
\begin{pmatrix}
1  & 0  \\ 
0  & 1  \\ 
\end{pmatrix} .
\end{equation} 
It is easy to see from the definition that 
$$
\det 
\begin{pmatrix} 
p_{\alpha , n-1}  & p_{\alpha , n } \\ 
q_{\alpha , n-1}  & q_{\alpha , n}  \\ 
\end{pmatrix}
\, = \, 
(-1)^{\sharp \{ 0 \le i \le n-1, \, \varepsilon_{\alpha, i} = 1\}}
$$
for $n \ge 1$. 
We put
$x_{\alpha, n} = A_{\alpha}^{n} (x_{\alpha,0})$, $n \ge 0$.  It follows from definitions that 
\begin{eqnarray}
x  - \frac{p_{\alpha , n } }{q_{\alpha , n} } 
& = & 
\frac{ p_{\alpha , n} + \varepsilon_{\alpha,n} p_{\alpha , n-1 } x_{\alpha, n} } 
{q_{\alpha , n}  + \varepsilon_{\alpha,n} q_{\alpha , n-1 } x_{\alpha, n} }  
 -  
 \frac{p_{\alpha , n } }{q_{\alpha , n} }  \notag \\
  {} & = & 
  \frac{(p_{\alpha , n-1 }q_{\alpha , n } - q_{\alpha , n-1 }p_{\alpha , n})
\varepsilon_{\alpha,n} x_{\alpha, n} }
{q_{\alpha , n } (q_{\alpha , n} + q_{\alpha , n-1} x_{\alpha, n}) } \notag \\
{} & = & 
\det (A_{\alpha, 1} \cdots A_{\alpha, n+1}) \, 
\frac{ x_{\alpha, n}} 
{q_{\alpha , n } (q_{\alpha , n} + q_{\alpha , n-1} x_{\alpha, n}) } \label{1}
\end{eqnarray}
and 
\[
\sign\left(x - \frac{p_{\alpha , n } }{q_{\alpha , n} }\right) 
\, = \, 
\sign(q_{\alpha , n} x \, - \, p_{\alpha , n} )
\, = \, 
(-1)^{\sharp \{ 0 \le i \le n, \, \varepsilon_{\alpha, i} = 1\}} .
\]
We define
\begin{displaymath}
\beta_{\alpha,n}x_{\alpha,0} \cdots x_{\alpha,n},\quad \beta_{\alpha,-1}=1 . 
\end{displaymath}

\begin{lem} \label{beta} 
For all $n \geq 1$,
\[
\beta_{\alpha,n} \, = \, 
|q_{\alpha, n} x \, - \, p_{\alpha, n}| 
\]
\end{lem}
\begin{proof} 
From
\begin{displaymath}
x=\frac{p_{\alpha,n}+p_{\alpha,n-1}\varepsilon_{\alpha,n}x_{\alpha,n}}{q_{\alpha,n}+q_{\alpha, n-1} \varepsilon_{\alpha,n} x_{\alpha,n}}  
\end{displaymath}
we get 
\begin{displaymath}
x_{\alpha,n} =-\varepsilon_{\alpha,n}\frac{p_{\alpha,n}-q_{\alpha,n} x}{p_{\alpha,n-1}-q_{\alpha, n-1} x} . 
\end{displaymath}
Then clearly $x_{\alpha,-1} x_{\alpha,0} \ldots x_{\alpha,n}=\abs{p_{\alpha,n}-q_{\alpha,n} x}$.
\end{proof}

From equation (\ref{1}), it follows that
\begin{displaymath} 
\beta_{\alpha,n}=\frac{x_{\alpha,n}}{q_{\alpha,n}+q_{\alpha,n-1}\varepsilon_{\alpha,n} x_{\alpha,n}} \; \Rightarrow \; \beta_{\alpha,n}=\frac{\beta_{\alpha,n+1}}{x_{\alpha,n+1}}=\frac{1}{q_{\alpha,n+1}+q_{\alpha,n}\varepsilon_{\alpha,n+1}x_{\alpha,n+1}} .
\end{displaymath}
Knowing that $x_{\alpha,n+1} \in (0,\bar{\alpha})$, we obtain
\begin{equation} \label{beta_q}
\frac{1}{1+\alpha} < \beta_{\alpha,n} q_{\alpha,n+1} < \frac{1}{\alpha} . 
\end{equation}

The proof of the following Lemma can be found in \cite{MCM}: 
\begin{lem} \label{rho}
Let $\bar{\alpha}=\max(\alpha,1-\alpha)$, $g=\frac{\sqrt{5}-1}{2}$, $\gamma=\sqrt{2}-1$. Then $\forall n \geq 0$,
\begin{displaymath}
\beta_{\alpha,n} \leq \bar{\alpha} \rho_{\alpha}^n, \qquad \frac{1}{q_{\alpha,n+1}}<(1+\alpha)\bar{\alpha} \rho_{\alpha}^n  
\end{displaymath}
where
\begin{itemize}
\item $\rho_{\alpha}=g$, \quad for $g < \alpha \leq 1$
\item $\rho_{\alpha}=\gamma$, \quad for $\gamma \leq \alpha \leq g$,  
\item $\rho_{\alpha}=\sqrt{1-2\alpha}$, \quad for $0 < \alpha < \gamma$. 
\end{itemize}
\end{lem}

\begin{rem}
The previous estimate is optimal for $\alpha \geq \gamma$ \cite{MCM}; however, we do not possess an optimal bound for $\alpha <\gamma$.
\end{rem}

\subsection{The $(\alpha, u)$-Brjuno functions}
It has already been remarked \cite{MMY2} that it is possible to extend the notion of Brjuno functions by replacing the logarithm with another function exhibiting a similar behaviour near $0$. 
Let $\alpha \in [0,1]$, and $u: (0,1) \to \R^+$ be a $C^1$ function such that
\begin{equation} \label{conditions}
\begin{split}
&\lim_{x \to 0^+} u(x)= \infty,\\
&\lim_{x \to 0^+} xu(x)< \infty,\\
&\lim_{x 	\to 0^+} x^2u'(x) < \infty.
\end{split}
\end{equation}
We define the generalized Brjuno function related to $\alpha$ and $u$: for $x \in {\mathbb I}_{\alpha}$,
\begin{displaymath}
B_{\alpha,u}(x)\doteqdot\sum_{n=0}^{\infty} \beta_{\alpha, n-1} u(x_{\alpha,n}) . 
\end{displaymath}
We extend $B_{\alpha,u}$ to the whole real line by setting 
\begin{align*}
& B_{\alpha,u}(x+1)=B_{\alpha,u}(x) \quad \forall x \in \R,\\
& B_{\alpha,u}(-x)=B_{\alpha,u}(x) \quad \quad \forall x \in (0, \min(\alpha,1-\alpha)) . 
\end{align*}
When $u(x)=-\log(x)$, we obtain the \emph{$\alpha$-Brjuno functions} \cite{MMY}:
\begin{displaymath}
B_{\alpha} (x) =\sum_{n=0}^{\infty}
 \beta_{\alpha,n-1} \log x_{\alpha, n} . 
\end{displaymath}

An irrational number $x$ will be called an \emph{$(\alpha,u)$-Brjuno number} if $B_{\alpha,u}(x_0) < \infty$, where $x_0=\abs{x-\floor{x+1-\alpha}}$. In this case we will also say that $x$ satisfies the \emph{$(\alpha,u)$-Brjuno condition}. 

\begin{rem} The first assumption in (\ref{conditions}) guarantees that the $(\alpha,u)$-Brjuno condition is not verified by all irrational numbers. The second is necessary in order to ensure that the set of $(\alpha,u)$-Brjuno numbers is nonempty. The third assumption is technical. 
\end{rem}

\begin{rem}
We remark here that in the case of the $\alpha$-Brjuno functions, we can replace $\sum_{n=0}^{\infty} \frac{\log(a_{\alpha,n+1})}{q_{\alpha,n}}$ with $\sum_{n=0}^{\infty} \frac{\log(q_{\alpha,n+1})}{q_{\alpha,n}}$ since the two series have the same set of convergence and their difference is bounded. In fact, the series $\sum\frac{\log(q_{\alpha,n})}{q_{\alpha,n}}$ is uniformly bounded for all $\alpha \in (0,1]$, see \cite{MMY}.
However, it is \emph{not} possible in general to replace $u(a_{\alpha, n+1}^{-1})$ with $u(q_{\alpha, n+1}^{-1})$, since the series $\sum\frac{u(q_{\alpha,n}^{-1})}{q_{\alpha,n}}$ might not be bounded.  
\end{rem}

\begin{figure}[htb] \label{sqrt_brjuno}
\begin{center}
\includegraphics[width=0.7\textwidth]{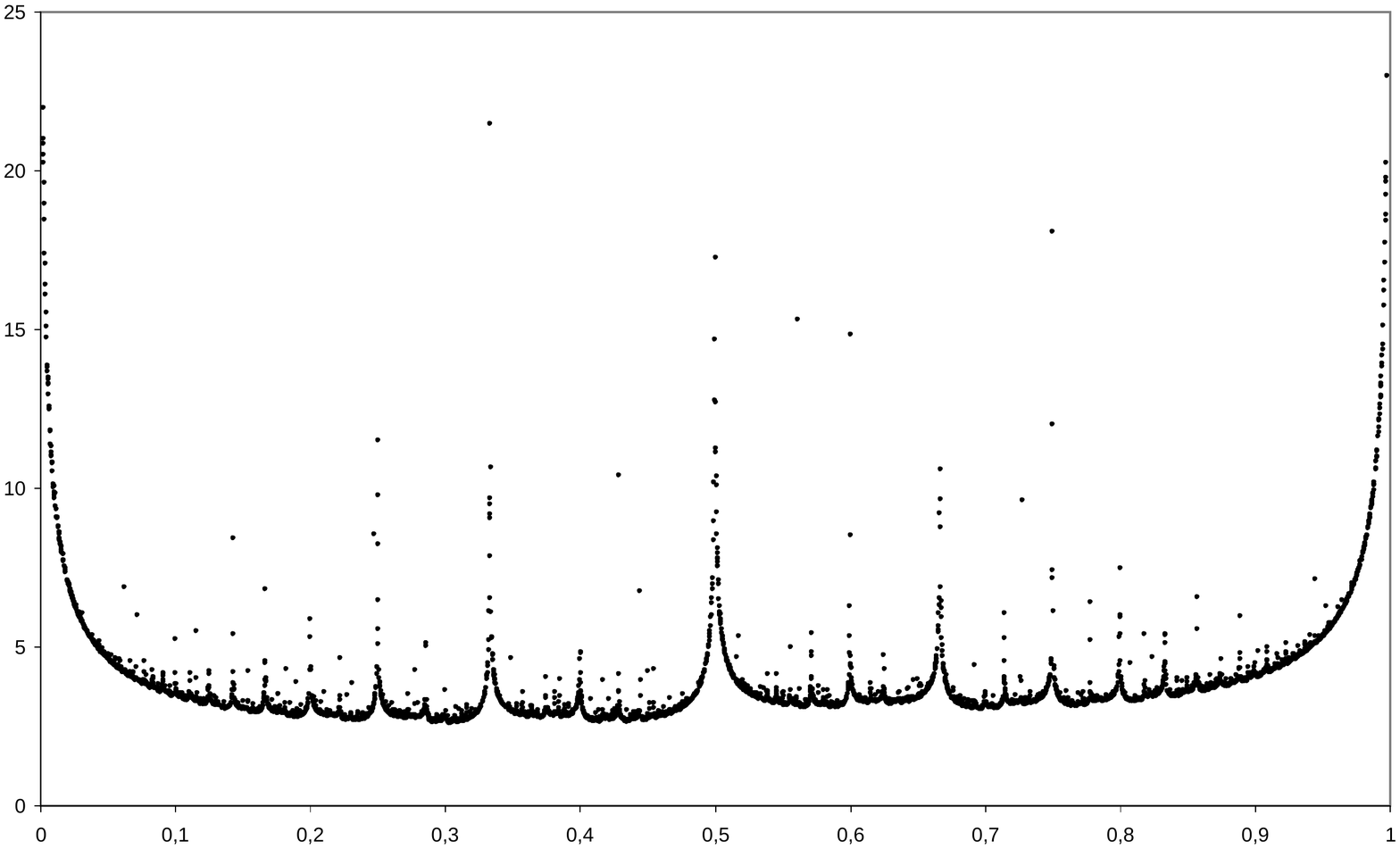}%
\caption{A computer simulation of the graph of $B_{(1,u)}$ when $u(x)=\frac{1}{\sqrt{x}}$.}
\end{center}
\end{figure}

Let $B_{\alpha,u}^N(x)=\sum_{n=0}^{\infty}
\beta_{\alpha,n}u(x_{\alpha,n})$ be the truncated $(\alpha,u)$-Brjuno series. 
\begin{prop} \label{q_series}
For any $\alpha, \, 0 < \alpha \le 1$, and for any $u$ satisfying the conditions (\ref{conditions}), there exists a positive 
constant $C_{\alpha, 1}$, depending on $u$, such that $\forall x \in \mathbb{I}_{\alpha}$, $\forall N>0$,
\[
\abs{  
B_{\alpha,u}^N(x) - \sum_{n=0}^{N} \frac{u(a_{\alpha, n+1}^{-1})}{q_{\alpha, n}} } < C_{\alpha, 1}  
\]
where $q_{\alpha,n}$ are the denominators of the $\alpha$-convergents of $x$.
\end{prop}

\begin{proof}[Proof of Proposition \ref{q_series}] 
We choose a positive constant $\delta (= \delta_{\alpha})$ sufficiently 
small. Let
\begin{align*}
& M_1(u)=\sup_{(\frac{\delta}{1+\delta},1)} \;u(x),\\
& M_2(u)=\sup_{(0,\delta)} \;xu(x),\\
& M_3(u)=\sup_{(0,\delta)} \;x^2 u'(x) . 
\end{align*}
We split the estimate in two cases: when $x_{\alpha,n}>\delta$, using Lemma \ref{rho} we get
\begin{displaymath}
\sum_{\substack{0 \leq n \leq N,\\ x_{\alpha,n}>\delta}} \beta_{\alpha, n-1} u(x_{\alpha,n}) < M_1\left(\sum_{n=1}^\infty \bar{\alpha} \rho_{\alpha}^{n-1}+1\right)=M_1 \left(\frac{\bar{\alpha}}{1-\rho_{\alpha}}+1\right)
\end{displaymath}
Moreover, observing that
\begin{displaymath}
a_{\alpha, n+1} = \left[ \frac{1}{x_{\alpha, n}} \right]_{\alpha}>\frac{1}{x_{\alpha,n}}-1 \quad \Rightarrow \quad \frac{1}{a_{\alpha,n+1}}>\frac{x_{\alpha,n}}{1+x_{\alpha,n}}>\frac{\delta}{1+\delta}
\end{displaymath}
we find, from Lemma \ref{rho}:
\begin{displaymath}
\sum_{\substack{0 \leq n \leq N,\\ x_{\alpha,n}>\delta}} \frac{u(a_{\alpha,n+1}^{-1})}{q_{\alpha,n}}\leq M_1(u)\left(1+\sum_{n=1}^{\infty} (1+\alpha) \bar{\alpha} \rho_{\alpha}^{n-1}\right)=M_1(u)\left(1+\frac{(1+\alpha)\bar{\alpha}}{1-\rho_{\alpha}}\right) 
\end{displaymath}
Suppose $x_{\alpha, n} < \delta$, which implies that $a_{\alpha, n+1}$ 
is not small, indeed $a_{\alpha, n+1} \geq 3$ is enough. We have
\begin{multline} \label{3}
\abs{\beta_{\alpha,n-1}u(x_{\alpha,n})-\frac{u(a_{\alpha,n+1}^{-1})}{q_{\alpha,n}}}< \\ <\beta_{\alpha,n-1}\abs{u(x_{\alpha,n})-u(a_{\alpha,n+1}^{-1})}+\abs{\beta_{\alpha,n-1}-\frac{1}{q_{\alpha,n}}}u(a_{\alpha,n+1}^{-1}) . 
\end{multline}
From Lagrange's Theorem we get 
\begin{displaymath}
\abs{u(x_{\alpha,n})-u(a_{\alpha,n+1}^{-1})}<\abs{u'(\xi)}\abs{x_{\alpha,n}-\frac{1}{a_{\alpha,n+1}}}
\end{displaymath}
for some $\xi$ between $x_{\alpha,n}$ and $a_{\alpha,n+1}^{-1}$. Since 
\begin{displaymath}
x_{\alpha,n},\frac{1}{a_{\alpha,n+1}},\xi \in \left(\frac{1}{a_{\alpha,n+1}+\alpha},\frac{1}{a_{\alpha,n+1}-1+\alpha}\right]
\end{displaymath}
we have
\begin{displaymath}
\abs{x_{\alpha,n}-\frac{1}{a_{\alpha,n+1}}}<\frac{1}{(a_{\alpha,n+1}-1+\alpha)(a_{\alpha,n+1}+\alpha)}<\frac{\xi^2}{1-\xi}<2\xi^2
\end{displaymath}
and
\begin{displaymath}
\sum_{\substack{0 \leq n \leq N \\ x_{\alpha,n}<\delta}}\beta_{\alpha,n-1}\abs{u(x_{\alpha,n})-u(a_{\alpha,n+1}^{-1})}<2M_3(u)\left(\frac{\bar{\alpha}}{1-\rho_{\alpha}}+1\right) . 
\end{displaymath}
For the second term in (\ref{3}), we observe that
\begin{align*}
&\abs{\beta_{\alpha,n-1}-\frac{1}{q_{\alpha,n}}}=\abs{\frac{1}{q_{\alpha,n}+q_{\alpha,n-1}\varepsilon_{\alpha,n}x_{\alpha,n}}-\frac{1}{q_{\alpha,n}}}=\frac{q_{\alpha,n-1}x_{\alpha,n}}{q_{\alpha,n}(q_{\alpha,n}+q_{\alpha,n-1}x_{\alpha,n})}=\\
&=\beta_{\alpha,n-1}\frac{q_{\alpha,n-1}x_{\alpha,n}}{q_{\alpha,n}}<\beta_{\alpha,n-1}x_{\alpha,n}
\end{align*}
and so 
\begin{multline*}
\abs{\beta_{\alpha,n-1}-\frac{1}{q_{\alpha,n}}}u(a_{\alpha,n+1}^{-1})<\\
<\beta_{\alpha,n-1}\left(x_{\alpha,n}u(x_{\alpha,n})
+x_{\alpha,n}\abs{u(x_{\alpha,n})-u(a_{\alpha,n+1}^{-1})}\right)<\\
<\beta_{\alpha,n-1}\left(M_2(u)+2M_3(u)\right) .
\end{multline*}
Finally,
\begin{displaymath} 
\sum_{\substack{0 \leq n \leq N \\ x_{\alpha,n}<\delta}}\abs{\beta_{\alpha,n-1}-\frac{1}{q_{\alpha,n}}}u(a_{\alpha,n+1}^{-1})<(M_2(u)+2M_3(u))\left(\frac{\bar{\alpha}}{1-\rho_{\alpha}}+1\right) . \qedhere
\end{displaymath}
\end{proof}

\begin{rem} If $B_{\alpha,u}(x) < \infty$, then 
$\sum_{n=0}^{\infty} \frac{u(a_{\alpha, n+1}^{-1})}{q_{\alpha, n}}$ is away from 
$B_{\alpha,u}(x)$ at most $C_{\alpha, 1}$. 
\end{rem}

\begin{theo} \label{bounded_difference}
Let $0 < \alpha \leq 1$ and assume that $u$ satisfies the conditions (\ref{conditions}). Then there
exists a positive constant $C_{\alpha,2}$, depending on $u$, such that 
\[
| B_{\alpha,u}(x) \, - \, B_{1,u}(x) | \, < \, C_{\alpha, 2} 
\]
whenever $B_{1,u}(x) < \infty$. 
\end{theo}

\begin{proof}
It has been proved in \cite{Nat2} that for any $\alpha \in (0,1]$, the $\alpha$-continued fraction expansion admits a \emph{Legendre constant}, that is, there exists a constant 
$k_{\alpha} > 0$ such that for any $n \ge 1$, 
\begin{equation} \label{b}
\left| x \, - \, \frac{p}{q} \right| 
\, < \, k_{\alpha} \frac{1}{q^{2}} \quad
\Rightarrow \quad \exists n \; :\; 
\frac{p}{q}
\, = \, 
\frac{p_{\alpha, n}}{q_{\alpha, n}} . 
\end{equation}
In particular, if $\frac{p}{q}=\frac{p_{1,l}}{q_{1,l}}$, there exists $n=n(x,l)$ such that $\frac{p_{1,l}}{q_{1,l}}=\frac{p_{\alpha,n}}{q_{\alpha,n}}$.
On the other hand, it is well-known that the Legendre constant when $\alpha=1$ is $\frac{1}{2}$; therefore, if for some $n\geq 0$ we have
\begin{equation} \label{a}
\left| x \, - \, \frac{p_{\alpha, n}}{q_{\alpha, n}} \right| 
\, < \frac{1}{2q_{\alpha, n}^{2}},
\end{equation}
then there exists $l= l(x, n)$ such that $\frac{p_{\alpha, n}}{q_{\alpha, n}} = \frac{p_{1, l}}{q_{1, l}}$.

We suppose that condition (\ref{a}) holds (the case when condition (\ref{b}) holds for $\frac{p}{q}=\frac{p_{1,l}}{q_{1,l}}$ is similar.) Since equation (\ref{b}) is verified for all constants $k_{\alpha}$ smaller than the Legendre constant, we can suppose without any loss of generality that $k_{\alpha}$ is small, for example $k_{\alpha} < \frac{1}{10}$, which implies that $x_{\alpha, n}$ is also small:
\begin{eqnarray*}
\abs{q_{\alpha,n}x-p_{\alpha,n}}=\frac{x_{\alpha,n}}{q_{\alpha,n}+q_{\alpha,n-1}\varepsilon_{\alpha,n} x_{\alpha,n}}<\frac{k_{\alpha}}{q_{\alpha,n}} \\
\Rightarrow x_{\alpha,n}<\frac{k_\alpha(q_{\alpha,n}+q_{\alpha,n-1}\varepsilon_{\alpha,n} x_{\alpha,n})}{q_{\alpha,n}}<2k_{\alpha}<\frac{1}{5} .
\end{eqnarray*}
In this case $\beta_{\alpha,n}=| q_{\alpha, n} x \, - \, p_{\alpha, n} |
\, =\, 
| q_{1, l}x \, - \, p_{1, l} |=\beta_{1,l}$; moreover, 
\begin{displaymath}
\beta_{\alpha, n}=\frac{1}{q_{\alpha, n} \frac{1}{x_{\alpha, n}} + q_{\alpha, n-1}\varepsilon_{\alpha,n} }, \quad \beta_{1,l}=\frac{1}{q_{1, l}  \frac{1}{x_{1, l} } + q_{1, l-1} }
\end{displaymath}
Hence we have 
\[
 q_{\alpha, n} \cdot \frac{1}{x_{\alpha, n} } + q_{\alpha, n-1} \varepsilon_{\alpha,n} 
\, = \, 
q_{1, l} \cdot \frac{1}{x_{1, l} } + q_{1, l-1} . 
\]
and so, dividing by $q_{1,l}=q_{\alpha,n}$, we get 
\[
\frac{1}{x_{\alpha,n}}-\frac{q_{\alpha,n-1}}{q_{\alpha,n}}\leq\frac{q_{1, l-1}}{q_{1, l}}  +  \frac{1}{x_{1, l}} \leq \frac{1}{x_{\alpha,n}}+\frac{q_{\alpha,n-1}}{q_{\alpha,n}}
\]
which implies 
\[
\left| \frac{1}{x_{\alpha, n}}  -  \frac{1}{x_{1, l}}\right|\leq \frac{q_{\alpha,n-1}}{q_{\alpha,n}}+\frac{q_{1,l-1}}{q_{1,l}}  
 <  2 .
\]
Then $\abs{x_{1,l}-x_{\alpha,n}}<2x_{1,l}x_{\alpha,n}<2x_{\alpha,n}<\frac{2}{5}$ and $x_{1,l}<3x_{\alpha,n}$, $x_{\alpha,n}<3x_{1,l}$.  
Suppose for example that $x_{\alpha,n}\leq x_{1,l}$ (the other case is similar); we want to estimate
\begin{multline*} 
\abs{\beta_{\alpha,n-1}u(x_{\alpha,n})-\beta_{1,l-1}u(x_{1,l})} \leq \\ \leq \beta_{\alpha,n-1}\abs{u(x_{\alpha,n})-u(x_{1,l})}+\abs{\beta_{\alpha,n-1}-\beta_{1,l-1}}u(x_{1,l}) .
\end{multline*} 
Since we are assuming that $x_{\alpha, n}$ is small, we have 
\begin{displaymath}
\beta_{\alpha,n-1}\abs{u(x_{\alpha,n})-u(x_{1,l})} \leq  \beta_{\alpha,n-1} u'(\xi) \abs{x_{1,l}-x_{\alpha,n}} \leq 2 \beta_{\alpha,n-1} u'(\xi)x_{1,l}x_{\alpha,n}
\end{displaymath}
with $\xi \in (x_{\alpha,n}, x_{1,l})$. Then 
\begin{displaymath}
\beta_{\alpha,n-1}\abs{u(x_{\alpha,n})-u(x_{1,l})} \leq 6 \beta_{\alpha,n-1}u'(\xi)\xi^2 \leq 6 \beta_{\alpha,n-1}M_3(u) . 
\end{displaymath}
On the other hand, we see that 
\begin{align*}
&  \abs{\beta_{\alpha,n-1}-\beta_{1,l-1}}u(x_{1,l})= \abs{\frac{1}{q_{\alpha,n}+q_{\alpha,n-1}\varepsilon_{\alpha,n}x_{\alpha,n}}-\frac{1}{q_{1,l}+q_{1,l-1}x_{1,l}}}u(x_{1,l})=\\
&=\frac{1}{q_{1,l}}\abs{\frac{1}{1+\frac{q_{\alpha,n-1}}{q_{\alpha,n}}\varepsilon_{\alpha,n} x_{\alpha,n}}-\frac{1}{1+\frac{q_{1,l-1}}{q_{1,l}}x_{1,l}}}u(x_{1,l})\leq \\
&\leq \frac{1}{q_{1,l}^2}\frac{q_{1,l-1}x_{1,l}+x_{\alpha,n}q_{\alpha,n-1}}{\left(1+\frac{q_{\alpha,n-1}}{q_{\alpha,n}}\varepsilon_{\alpha,n}x_{\alpha,n}\right)\left(1+\frac{q_{1,l-1}}{q_{1,l}}x_{1,l}\right)}u(x_{1,l})\leq \frac{x_{1,l}+\abs{x_{\alpha,n}}}{q_{1,l}}u(x_{1,l})<\\
&<\frac{2x_{1,l}}{q_{1,l}}u(x_{1,l})\leq 2M_2(u)(1+\alpha) \bar{\alpha}\rho_{\alpha}^{l-1},
\end{align*}
where the last equality holds for $l \geq 1$. 
\bigskip \par
If neither (\ref{a}) nor (\ref{b}) hold, then
\begin{align*}
&\abs{x-\frac{p_{\alpha,n}}{q_{\alpha,n}}}=\frac{x_{\alpha,n}}{q_{\alpha,n}(q_{\alpha,n}+q_{\alpha,n-1}\varepsilon_{\alpha,n}x_{\alpha,n})}\geq \frac{k_{\alpha}}{q_{\alpha,n}^2} \\
&\Rightarrow x_{\alpha,n}\geq k_{\alpha} \left(1-\frac{q_{\alpha,n-1}}{q_{\alpha,n}}x_{\alpha,n}\right) \geq k_{\alpha} (1-x_{\alpha,n}) 
\end{align*}
and so $x_{\alpha,n}\geq \frac{k_{\alpha}}{1+k_{\alpha}}$ is bounded away from $0$; the same is true for $x_{1,l}$. In this case for $l, n \geq 1$ we have
\begin{displaymath}
\abs{\beta_{\alpha,n-1}u(x_{\alpha,n})-\beta_{1,l-1}u(x_{1,l})} \leq M_4(u) (\bar{\alpha}\rho_{\alpha}^{n-1}+ g^{l-1}),
\end{displaymath}
where 
\begin{displaymath}
M_4(u)=\max\limits_{\left[\frac{k_{\alpha}}{1+k_{\alpha}},1\right]} u(x) . \qedhere
\end{displaymath}
\end{proof}

\begin{rem} From Proposition \ref{q_series} and Theorem \ref{bounded_difference}, we see that, similarly to what happens in the case $u(x)=-\log x$, the set
$\{x \in \R \, : \, B_{\alpha,u}(x) < \infty \}$ coincides with the set of $(1,u)$-Brjuno 
numbers for any $\alpha \in (0, 1]$.
\end{rem}

It is known \cite{BDV} \cite{LM} that the maps $A_{\alpha}$ admit a unique absolutely continuous invariant density $\frac{d\mu_{\alpha}(x)}{dx}$ which is bounded from above and from below by constants depending on $\alpha$.\footnote{The density is known explicitly when $\alpha \in [\gamma,1]$, see \cite{Nak} and \cite{MCM}.}\\
Following \cite{MCM}, we consider the linear space
\begin{multline*}
X_{\alpha}=\left\{f: \R \to \R \text{ measurable }\;|\; f(x+1)=f(x) \; \forall x,\right. \\ \left. \; f(x)=f(-x) \text{ for } x \in (0,\min(\alpha,1-\alpha))\right\} . 
\end{multline*}
We call $X_{\alpha,p}$ the space $X_{\alpha}$ endowed with the $L^p$ norm
\begin{displaymath}
\norm{f}_p=\left(\int_0^{\bar{\alpha}}\abs{f}^p d\mu_{\alpha}(x)\right)^{\frac{1}{p}} . 
\end{displaymath}
For $f \in X_{\alpha,p}$, and for $x \in (0, \bar{\alpha})$, we define 
\begin{displaymath}
(T_{\alpha}f)(x)=xf (A_{\alpha}(x))
\end{displaymath}
and we extend the domain of $T_{\alpha}f$ with the same periodicity and parity conditions in the definition of $X_{\alpha}$. Then $T_{\alpha}f \in X_{\alpha,p}$. 
It is easy to check that for $\alpha \in \left[\frac{1}{2},1\right]$, for all $x \in(0,\bar{\alpha})$ one has
\begin{displaymath}
(1-T_{\alpha})B_{\alpha,u}(x)=u_{\alpha}(x),
\end{displaymath}
where $u_{\alpha}$ is the extension of $u$ to the whole real line which coincides with $u$ on $(0,\bar{\alpha})$ and belongs to $X_{\alpha}$. Similarly, for all $x \in (0,1-\alpha)$ and for $\alpha \in \left(0,\frac{1}{2}\right)$ one has
\begin{displaymath}
(1-T_{\alpha} I)B_{\alpha,u}(x)=u_{\alpha}(x), 
\end{displaymath} 
where $I$ denotes the involution $(If)(x)=f(-x)$. 
The spectral radius of $T_{\alpha}$ in $X_{\alpha,p}$ is bounded by $\rho_{\alpha}$, and so $1-T_{\alpha}$ is invertible \cite{MCM}. As a consequence, we have the following
\begin{prop}
For $\alpha \in (0,1]$, if $u_{\alpha} \in L^p(0,\bar{\alpha})$, then $B_{\alpha,u} \in L^p(\mathbb{T}^1)$.
\end{prop}
\begin{rem}
Since $x \in C^1(0,\bar{\alpha})$ the assumption $u_{\alpha} \in L^p$ is really a growth assumption as $x \to 0$. If $u(x)=x^{-\frac{1}{\sigma}}$, $\sigma>1$ then $u_{\alpha} \in L^p$ for all $p<\sigma$. In this case the set of $(\alpha,u)$-Brjuno numbers includes the set of diophantine numbers with exponent $\tau<\sigma$ and is included in the set of diophantine numbers with exponent $\tau=\sigma$ \cite{MMY2}.
\end{rem}

\section{The semi-Brjuno function}
\subsection{Some definitions}
In the sequel we focus on the relations between the standard Brjuno function $B_1$ and the Brjuno function $B_0$ associated to the by excess continued fraction approximation $A_0(x): (0,1] \to (0,1]$, $A_0(x)=\floor{\frac{1}{x}+1}-\frac{1}{x}$.\\
We recall some basic properties of these expansions: let $x^n = A_0^n (x)$ for $n \ge 0$, $x_{-1} = 1$, and let
\[
b_n \, = \floor{\frac{1}{x_{n-1}}+1} \qquad \mbox{for} \,\, n \ge 1. 
\]
be the elements of the \vv{--} continued fraction expansion of $x$:
\begin{equation} \label{-expansion}
x \, = \,   \confrac{1}{b_1} - \confrac{1}{b_2} -
\confrac{1}{b_3} - \cdots \, .
\end{equation}
For simplicity of notation, we denote by $\frac{p_n^*}{q_n^*}=\frac{p_{0,n}}{q_{0,n}}$ the $0$-convergents of $x$ defined in equations (\ref{convergents}) and (\ref{p0}), so that
\[
\frac{p_{n}^{\ast}}{q_{n}^{\ast}}  =  
\confrac{1}{b_1} - \confrac{1}{b_2} -
\confrac{1}{b_3} - \cdots - \confrac{1}{b_n}\,
\]
and $p_{n}^{\ast}q_{n-1}^{\ast}  - p_{n-1}^{\ast}q_{n}^{\ast} =  1$. Moreover we have 
\begin{displaymath}
 x = \frac{p_{n}^{\ast} - p_{n-1}^{\ast} x_n}
                {q_{n}^{\ast} - q_{n-1}^{\ast} x_n}, \quad \quad
x  -  \frac{p_{n}^{\ast}}{q_{n}^{\ast}} =  
\frac{x_n}{q_{n}^{\ast}(q_{n}^{\ast} - q_{n-1}^{\ast} x_n)} >  0,
\end{displaymath}
that is, the approximation by $A_0$ is \emph{one-sided}. Let 
\begin{equation} \label{eq1}
\beta_n^*=q_{n}^{\ast} x  -  p_{n}^{\ast}  =  x_{0} x_{1} \cdots x_{n}
\end{equation}
for $n \ge 1$.  From the definition, it follows that $q_{n}^{\ast} >  q_{n-1}^{\ast}$
for $n \ge 0$ and $q_{n}^{\ast} > 2 q_{n-1}^{\ast}$ if $b_n \ge 3$. It is easy to see that $b_1 = k$ if and only if 
$x \in \left(\frac{1}{k}, \, \frac{1}{k-1}\right]$.  If $x_{n-1} \in 
\left(1 - \frac{1}{l}, 1 - \frac{1}{l+1} \right]$ for $l \ge 2$, then 
$b_n = 2$ and $ x \in \left(1 - \frac{1}{l-1}, \, 1 - \frac{1}{l}\right]$. 
Thus, 
\[
b_n = b_{n+1} = \cdots = b_{n+l-2} = 2 
\]
when $x_{n-1} \in \left(1 - \frac{1}{l},  1 - \frac{1}{l+1} \right]$ for 
$l \ge 2$.  In this case, 
\[
q_{n+k}^{\ast} \, = \, q_{n-1}^{\ast} \, + \, 
                          (k+1)(q_{n-1}^{\ast} - q_{n-2}^{\ast}) 
\]
for $0 \le k \le l - 2$, and we may have 
\[
\sum_{n=1}^{\infty} \frac{1}{q_{n}^{\ast}} \, = \, \infty
\]
when we have ``long" consecutive sequences of ``$2$" infinitely often 
in the sequence $(b_n)$, as opposed to the case of the $\alpha$-continued fraction 
expansions with $\alpha \in (0,1]$, for which we always have 
\[
\sum_{n=1}^{\infty} \frac{1}{q_{\alpha, n}} \, < \, \infty
\]
for any irrational number $x$.  

As we will see, a precise \vv{dictionary} between the coefficients $\{a_n \, : \, n \ge 1 \}$ of the regular continued fraction expansion 
of $x$ and the by excess coefficients $\{b_n \, : \, n \ge 1 \}$ is available 
(see \cite{KN}). 
To begin with, we recall how to obtain the $1$-expansion of $x$ given the $1$-expansion of $1-x$:
\begin{lem} \label{lemma1}
Suppose $x \in \left(0,\frac{1}{2}\right)$ irrational. Then 
\[
x = \confrac{1}{a_1} + \confrac{1}{a_2 } + \confrac{1}{a_3} +
\confrac{1}{a_4} + \cdots \, . 
\]
implies
\[
1 - x = \confrac{1}{1} + \confrac{1}{a_1 - 1} + \confrac{1}{a_2} +
\confrac{1}{a_3} + \cdots \, . 
\]
\end{lem}
\begin{lem} \label{lemma2}
Let $x \in \left(\frac{1}{2},1\right)$ whose $0$-expansion is of the form (\ref{-expansion}), and $\{n_i\}_{i\geq 1}$ a sequence such that
\begin{equation} \label{eq2}
\begin{array}{l}
b_1 = b_2 = \cdots = b_{n_1 - 1} = 2 \\
b_{n_1} > 2 \\
b_{n_1 + 1} = \cdots = b_{n_1 + n_2 - 1} = 2 \\
b_{n_1 + n_2 } > 2 \\
b_{n_1 + n_2 + 1} = \cdots = b_{n_1 + n_2 + n_3 - 1} = 2 \\
\cdots 
\end{array}
\end{equation}
($n_i = 1$ is possible for $i \ge 2$).  Then the $1$-expansion of $x$ is 
\[
1 - x \, = \, 
\confrac{1}{n_1} + \confrac{1}{b_{n_1}-2} + \confrac{1}{n_2} + 
\confrac{1}{b_{n_1 + n_2}-2} + \confrac{1}{n_3} + \cdots \, . 
\]
\end{lem}
Lemma \ref{lemma1} and \ref{lemma2} imply the following \cite{KN} 
\begin{prop} \label{prop1}
Suppose $x \in \left(\frac{1}{2}, 1\right)$.  If (\ref{eq2}) holds, then 
\[
x \, = \, 
\confrac{1}{1} + \confrac{1}{n_1 - 1} + \confrac{1}{b_{n_1}-2} + 
\confrac{1}{n_2} + 
\confrac{1}{b_{n_1 + n_2}-2} + \confrac{1}{n_3} + \cdots \, . 
\]
\end{prop}

In the case $x \in \left(0, \frac{1}{2}\right)$, then $n_1 = 1$ in (\ref{eq2}) and a similar result holds:
\begin{prop} \label{prop2}
If $x \in \left(0, \frac{1}{2}\right)$, then 
\[
x \, = \, 
 \confrac{1}{b_{n_1} - 1} + \confrac{1}{n_2} + 
\confrac{1}{b_{n_1 + n_2}-2} + \confrac{1}{n_3} + \cdots \, . 
\]
\end{prop}

\subsection{Semi-Brjuno numbers}
An irrational number $x \in(0, 1)$ is said to be a \emph{semi-Brjuno number} if 
\[
\sum_{n=0}^{\infty} \frac{\log(b_{n+1}- 1)}{q_{n}^{\ast}} \, < \, \infty .
\]
For an irrational number $x$, $x$ is said to be a semi-Brjuno number when its 
fractional part $x - \floor{x}$ is a semi-Brjuno number.  
For $x \in (0, 1)$ irrational, we consider the $0$-Brjuno function\footnote{Remark that the finiteness of $B_0$ on the rationals depends on the choice of the value of $A_0$ on the negligible set $\left\{\frac{1}{n}\right\}_{n \geq 1}$. If instead of $A_0$ we consider another version of the by-excess map (see \cite{It}), that is
$$S(x)=\left \lceil \frac{1}{x}\right \rceil -\frac{1}{x},$$ then $S\left(\frac{1}{n}\right)=0$ for all $n \geq 1$, and the iterates of any rational number would be all zero after a certain index, so that $B_0$ would diverge on all the rationals. However, with our definition $A_0\left(\frac{1}{n}\right)=1$ for all $n \geq 1$, and $B_0$ is finite on $\Q$. We also note that under the first assumption the by-excess expansion of a rational number is finite, while under the second it terminates in an infinite sequence of \vv{2}.
} 
\[
B_0(x) = \sum_{n=0}^{\infty} \beta_{n-1}^* \log \frac{1}{x_n} 
\]
where the coefficients $\beta^*_{n}$ are defined as in (\ref{eq1}). Equivalently, we can write
\[
B_0(x) = \sum_{n=0}^{\infty} \beta_{n-1}^*
\log \frac{\beta_{n-1}^*}{\beta_n^*} . 
\]
We extend the definition to all $x \in \R \setminus \Q$ by putting $B_0(x) \doteqdot B_0(x - \floor{x})$.  
We call $B_0$ the \emph{semi-Brjuno function}.  It is easy to see that for $x \in (0,1)$ irrational, $B_0$ satisfies the functional equation
\[
B_0(x) \, = \, \log \frac{1}{x} \, + \, x B_0\left(-\frac{1}{x}\right) .
\]

\begin{figure}[htb] \label{B0.pdf}
\begin{center}
\includegraphics[width=0.7\textwidth]{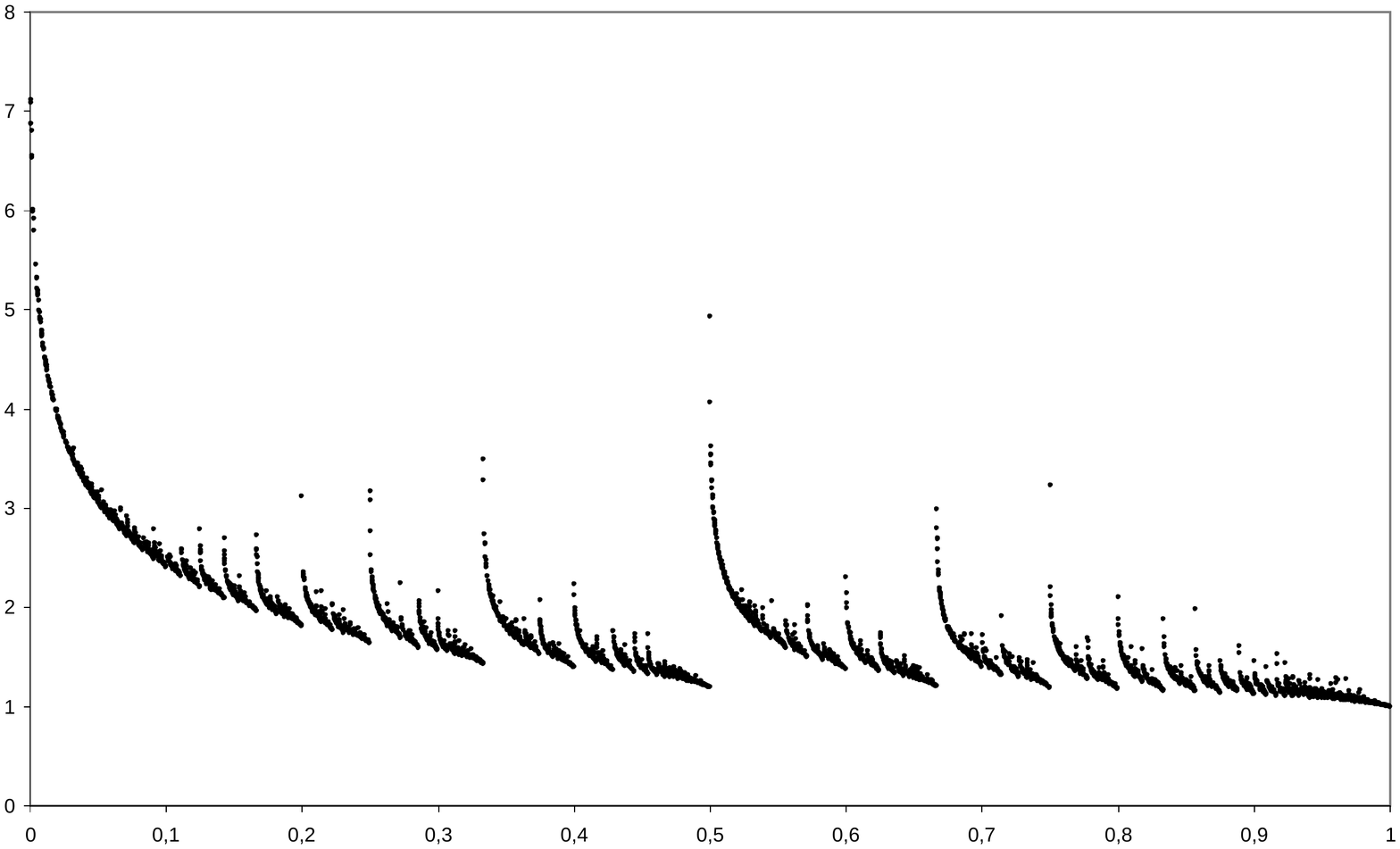}
\caption{A computer simulation showing the graph of $B_0$.}
\end{center}
\end{figure}

\begin{prop} \label{theorem1} 
An irrational number $x$ is a semi-Brjuno number if and 
only if $B_0(x) < \infty$.  
\end{prop}

Actually we will prove a stronger result:
\begin{theo} \label{theorem4}
There exists a fixed constant $C_0$ such that $\forall x \in (0,1)$,
$$ \abs{B_0(x)-\sum_{n=0}^\infty \frac{\log(b_{n+1}-1)}{q_n^*}}<C_0 .  $$
\end{theo}

\begin{proof}[Proof of Theorem \ref{theorem4}]
We split the series in the following way:
\[
B_0(x) \, = \, \sum_{n \in I_*(x)} \beta_{n-1}^* 
                                              \log \frac{1}{x_n} 
\, + \, \sum_{n \in I_{**}(x)} \beta_{n-1}^* 
                               \log \frac{1}{x_n},
\]
where
\[
\begin{array}{l}
I_*(x)=\left\{n \;|\; x_n \in \left(\frac{1}{2}, 1\right]\right\}=\left\{n \;|\; b_{n+1}(x)=2\right\}, \\
I_{**}(x)=\left\{n \;|\; x_n \in \left(0, \frac{1}{2}\right]\right\}=\left\{n \;|\; b_{n+1}(x)>2\right\} .
\end{array} 
\]
We also have
\[
\sum_{n=0}^{\infty} \frac{\log (b_{n+1} - 1)}{q_n^{\ast}} \, = \, 
\sum_{n \in I_*(x)}\frac{\log (b_{n+1} - 1)}{q_n^{\ast}}
\, + \, \sum_{n \in I_{**}(x)} \,\,\frac{\log (b_{n+1} - 1)}{q_n^{\ast}} \, .
\]
It is trivial that 
\[
\sum_{n \in I_{*}(x)} \frac{\log (b_{n+1} - 1)}{q_n^{\ast}} \, = \, 0 .
\]
We will show that $\sum_{n \in I_*(x)} \beta_{n-1}^* \log \frac{1}{x_n} < \infty$ is uniformly bounded by a fixed constant.
We consider $m \ge 1$ such that 
\[
x_m \leq \frac{1}{2}, \, x_{m+1} > \frac{1}{2}, \, x_{m+2} > \frac{1}{2}, \, 
\ldots \, , x_{m+l} > \frac{1}{2}, \, x_{m+l+1} \leq \frac{1}{2}. 
\]
In this case, \\
$x_{m+l} \in  \left(1 - \frac{1}{2}, 1 - \frac{1}{3}\right], 
 x_{m+l-1} \in \left(1 - \frac{1}{3}, 1 - \frac{1}{4}\right], 
 \cdots ,
 x_{m+1} \in \left(1 - \frac{1}{l+1}, 1 - \frac{1}{l+2}\right]$
For $0 < t \le l-1$, 
\begin{multline*} 
x_{-1} x_0 x_1 \cdots x_m x_{m+1} \cdots x_{m+t} \log \frac{1}{x_{m+t+1}} <\\
< (x_{-1} x_0 x_1 \cdots x_m ) \left(1 - \frac{1}{l+2}\right)\left(1 - \frac{1}{l+1}\right) 
            \cdots \left(1 - \frac{1}{l+3-t}\right) \frac{1}{l+1-t}= \\
= (x_{-1} x_0 x_1 \cdots x_m ) \frac{1}{l+2} \left(\frac{l+2-t}{l+1-t}\right)<(x_{-1} x_0 x_1 \cdots x_m ) \frac{1}{l+2};
\end{multline*}
for $t = 0$, 
\[
x_{-1} x_0 x_1 \cdots x_m \log \frac{1}{x_{m+1}} \, < \, 
x_{-1} x_0 x_1 \cdots x_m \frac{1}{l} .
\]
Thus we have 
\[
\sum_{t=0}^{l-1} x_{-1} x_0 x_1 \cdots x_{m+t} \log \frac{1}{x_{m+t+1}} 
\, < \, 2 (x_{-1} x_0 x_1 \cdots x_{m}). 
\]
We define 
\[
\left\{
\begin{array}{l}
s_1 \, = \, \min \{n \ge 0 \, : \, x_n \leq \frac{1}{2}, x_{n+1} > \frac{1}{2} \} 
\\
s_{t+1} \, = \, \min \{n > s_t \, : \, x_n \leq \frac{1}{2}, x_{n+1} > 
\frac{1}{2} \} , \,\,\, t \ge 1 .
\end{array} \right.
\]
Then 
\begin{multline*} 
\sum_{n \in I_*(x)} \,\, x_{-1} x_0 x_1 \cdots x_{n-1} \log\frac{1}{x_n} 
\, < \, 2 \sum_{j=1}^{\infty} x_{-1} x_0 x_1 \cdots x_{s_j}\leq \\
\leq 2 \sum_{j=1}^{\infty} x_{s_1} \cdots x_{s_j}< 2 \sum_{j=1}^\infty \frac{1}{2^j}=2 . 
\end{multline*}
Now we want to show that 
$$\abs{\sum_{n \in I_{**}(x)} \beta_{n-1}^* \log\frac{1}{x_n}-\sum_{n \in I_{**}(x)} \frac{\log(b_{n+1}-1)}{q_n*}}$$ is bounded by a constant which does not depend on $x$. Let 
\begin{equation} \label{t_j}
t_j = \sum_{i=1}^{j} n_i  - 1,
\end{equation}
so that $I_{**}(x)=\{t_1,t_2,t_3,\ldots\}$. Then $\forall k\geq 1$, $q_{t_{k+1}}^* \geq q_{t_k+1}^*=b_{t_k+1}q_{t_k}^*-q_{t_k-1}^* \geq 2 q_{t_k}^*$, and so $q_{t_k}^* \geq 2^k$, and
\begin{equation} \label{sqrt}
\sum_{n \in I_{**}(x)} \frac{1}{q_n^{\ast}}<2, \quad \quad \sum_{n \in I_{**}(x)} \frac{\log(q_n^*)}{q_n^{\ast}}<\sum_{n \in I_{**}(x)} \frac{1}{\sqrt{q_n^{\ast}}}<\frac{\sqrt{2}}{\sqrt{2}-1} . 
\end{equation}
Then, observing that  
\begin{displaymath}
\frac{\log q_n^*}{q_n*}+\frac{\log(b_{n+1}-1)}{q_n^*} \leq \frac{\log q_{n+1}^*}{q_n*} < \frac{\log q_n^*}{q_n*}+\frac{\log b_{n+1}}{q_n^*} 
\end{displaymath}
this is equivalent to showing that
\begin{equation} \label{MMY}
\abs{\sum_{n \in I_{**}(x)} \beta_{n-1}^* \log\frac{1}{x_n}-\sum_{n \in I_{**}(x)} \frac{\log(q_{n+1}^*)}{q_n*}}
\end{equation}
is uniformly bounded.
The proof is very similar to its analogue for the standard Brjuno function \cite{MMY}: from the relation 
$$q_n^* \beta_{n-1}^*-q_{n-1}^*\beta_n^*=1,$$
which can be checked easily, we get $\frac{1}{q_n^*}=\beta_{n-1}^*-\frac{q_{n-1}^* \beta_n^*}{\beta_n^*}$. Then expression (\ref{MMY}) becomes
\begin{align*}
&\abs{\sum_{n \in I_{**}(x)} \left(\beta_{n-1}^* \log \frac{\beta_n^*}{\beta_{n-1}^*}+ \log q_{n+1}^*\left(\beta_{n-1}^*-\frac{q_{n-1}^*}{q_n^*}\beta_n^*\right)\right)} \leq\\ 
& \leq\sum_{n \in I_{**}(x)} \abs{ \beta_{n-1}^* \log(\beta_n^* q_{n+1}^*)-\beta_{n-1}^* \log\beta_{n-1}^*-\frac{q_{n-1}^* \beta_n^*}{q_n*} \log(q_{n+1}^*)} . 
\end{align*}
For the first term, observe that when $b_{n+1}>2$, $x_n<\frac{1}{2}$ and 
\begin{align*}
&\beta_{n-1}^* = q_{n-1}^{\ast} x \, - \, p_{n-1}^{\ast}  = \frac{x_{n-1}}{q_{n-1}^{\ast} \, - \, q_{n-2}^{\ast} x_{n-1} }  = \frac{1}{q_{n}^{\ast} \, - \, q_{n-1}^{\ast} x_n }  < \frac{1}{\frac{1}{2} q_{n}^{\ast} }  = \frac{2}{q_{n}^{\ast}},\\
&\beta_n^*=\frac{1}{q_{n+1}^*-q_n^*x_{n+1}} \leq \frac{1}{q_{n+1}^*-q_n^*}\leq \frac{2}{q_n^*} \\
& \Rightarrow \sum_{n \in I_{**}(x)} \beta_{n-1}^* \log(\beta_n^* q_{n+1}^*) \leq \sum_{n \in I_{**}(x)} \frac{2\log 2}{q_n^*} \leq 4\log 2 . 
\end{align*}
From equation (\ref{sqrt}), it follows easily that $$\sum_{n \in I_{**}(x)} \beta_{n-1}^* \log \beta_{n-1}^* \leq \sum_{n \in I_{**}(x)}\frac{2\log(q_n^*/2)}{q_n^*} \leq \sqrt{2}\sum_{n \in I_{**}(x)} \frac{1}{\sqrt{q_n^*}}< \frac{2}{\sqrt{2}-1} . $$  
Finally, remarking that the function $\frac{\log(x)}{x}$ is decreasing for $x \geq e$, we have
\begin{align*}
& \sum_{n \in I_{**}(x)} \frac{q_{n-1}^* \beta_n^*}{q_n^*}\log(q_{n+1}^*) \leq \sum_{n \in I_{**}(x)} \beta_n^* \log q_{n+1}^* \leq \sum_{n \in I_{**}(x)} \frac{2 \log q_{n+1}^*}{q_{n+1}^*} \leq \\
& \leq \log 2+ \sum_{n \in I_{**}(x)} \frac{2 \log q_{n}^*}{q_{n}^*} \leq \log 2+\frac{2 \sqrt{2}}{\sqrt{2}-1}
\end{align*}
which completes the proof.
\end{proof}

Now we show the following 
\begin{prop} \label{theorem2} 
For any irrational number $x$, $x$ is a Brjuno number if and 
only if $x$ and $-x$ are semi-Brjuno numbers. 
\end{prop}
Again, we will actually prove a stronger result:
\begin{theo} \label{theorem5}
Let $B_0^+(x)=B_0(x)+B_0(-x)$ be the even part of $B_0$. 
Then the function $B_1(x)-B_0^+(x)$ is bounded.
\end{theo}

\begin{figure}[htb] \label{B1andB0even}
\begin{center}
\includegraphics[width=0.7\textwidth]{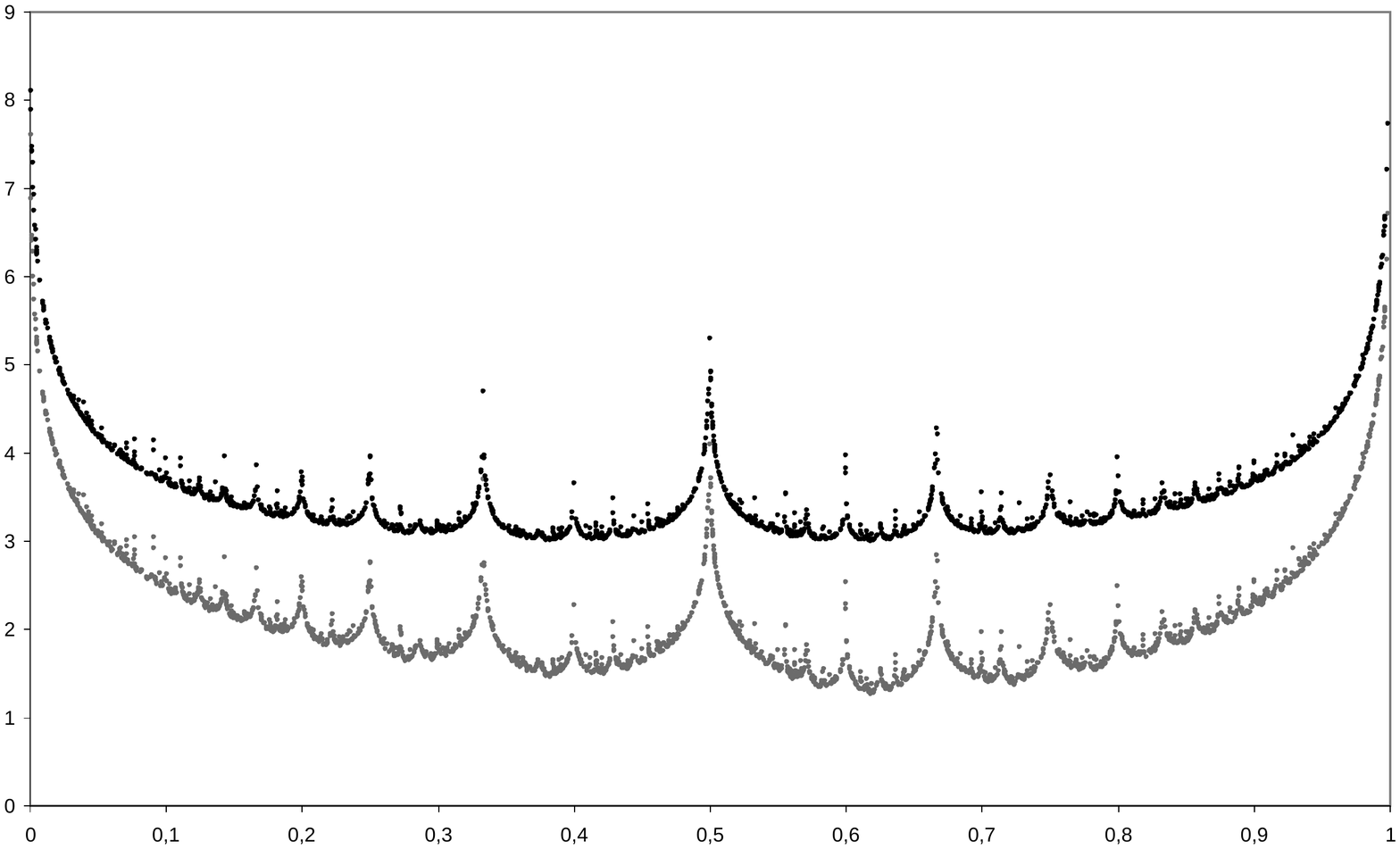}
\caption{The graphs of the functions $x \mapsto B_0(x)+B_0(-x)$ (black) and $x \mapsto B_1(x)$ (grey).}
\end{center}
\end{figure}

\begin{figure}[htb] \label{B1-B0even}
\begin{center}
\includegraphics[width=0.7\textwidth]{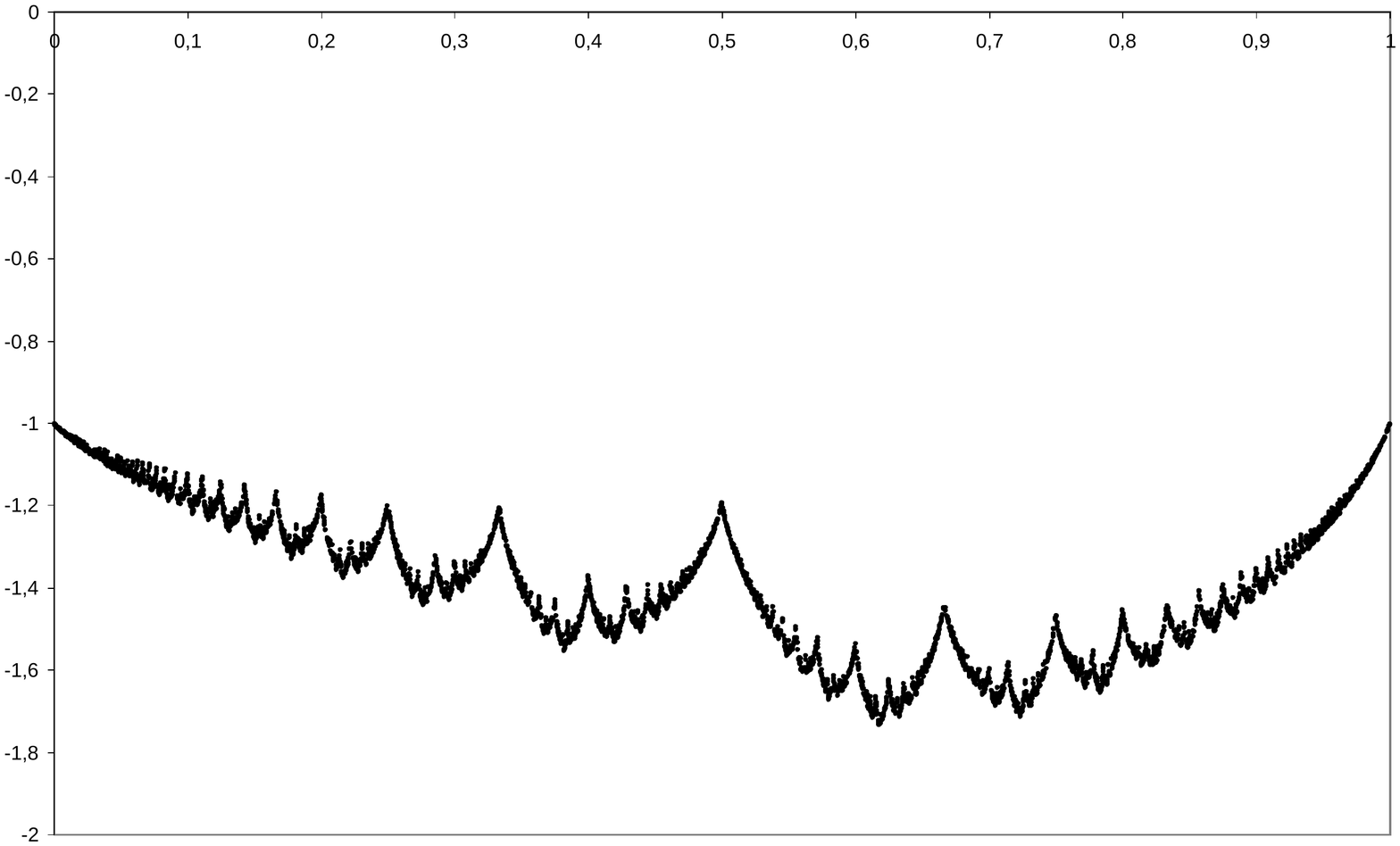}
\caption{The graph of the function $x \mapsto B_1(x)-(B_0(x)+B_0(-x))$.}
\end{center}
\end{figure}

For the proof of this theorem, we use the following 
\begin{lem} \label{lemma3}
There exists a constant $C_1$ such that for any irrational number $x \in (0, 1)$, 
\[
\abs{ \sum_{n=0}^{\infty} \frac{\log q_{n+1}}{q_n} -
\sum_{n=0}^{\infty} \frac{\log a_{n+1}}{q_n}} < C_1 .
\]
\end{lem}
\begin{proof}
Since $q_{n+1} = a_{n+1} 
q_{n} + q_{n-1}$, we see
\[
\sum_{n=0}^{\infty} \frac{ \log (a_{n+1} q_n + q_{n-1})}{q_n} 
\, < \, 
2 \left( \sum_{n=1}^{\infty} \frac{\log a_{n+1}}{q_n} + 
 \sum_{n=1}^{\infty} \frac{\log q_{n}}{q_n}
\right).
\]
It is known \cite{MMY} that for any $x$ irrational,
\begin{displaymath}
\sum_{n=1}^\infty \frac{\log q_n}{q_n} \leq \frac{2}{e}\left(3+\frac{\sqrt{2}g}{1-\sqrt{g}}\right) 
\end{displaymath}
where $g$ denotes as usual the small Golden ratio. 
\end{proof}

\begin{proof}[Proof of Theorem \ref{theorem5}] Observe that
\[
\sum_{n=0}^{\infty} \frac{\log a_{n+1}}{q_n} < \infty 
\quad \Leftrightarrow \quad \sum_{n=0}^{\infty} \frac{\log a_{2n+1}}{q_{2n}} < \infty
\,\, \mbox{ and } \,\, 
\sum_{n=0}^{\infty} \frac{\log a_{2n+2}}{q_{2n+1}} < \infty  .
\]
We will show that $\exists C_{2}>0$ such that $\forall x$ irrational, 
\[
\abs{\sum_{n=0}^{\infty} \frac{\log a_{2n+1}}{q_{2n}} -
\sum_{n=0}^{\infty} \frac{\log(b_{n+1}- 1)}{q_{n}^{\ast}}}<C_{2} .
\]
Suppose that $\frac{1}{2} < x < 1$.  Recall the definition of the sequence $\{t_j\}$ in (\ref{t_j}).
We have seen in the proof of Theorem \ref{theorem4} that 
\begin{displaymath}
\sum_{n=0}^{\infty} \frac{\log(b_{n+1}- 1)}{q_{n}^{\ast}} =
\sum_{j=1}^{\infty} \frac{\log(b_{t_j +1}- 1)}{q_{t_j}^{\ast}} .
\end{displaymath}
Moreover, because of the estimate (\ref{sqrt}), we have 
\begin{displaymath}
\abs{\sum_{j=1}^{\infty} \frac{\log(b_{t_j +1}- 1)}{q_{t_j}^{\ast}}-\sum_{b_{{t_j}+1} \geq 4} \frac{\log(b_{t_j +1})}{q_{t_j}^{\ast}}}=\sum_{b_{{t_j}+1} =3}\frac{\log 2}{q_{t_j}^*} \leq 2 \log 2 .
\end{displaymath}
When $b_{{t_j}+1} \geq 4$, since $x_{t_j} < \frac{1}{3}$, we get 
\begin{align*}
0 < x - \frac{p_{t_j}^{\ast}}{q_{t_j}^{\ast}} =\frac{x_{t_j}}{q_{t_j}^*(q_{t_j}^*-q_{t_j-1}^*x_{t_j})} < \frac{1}{3q_{t_j}^*}\frac{1}{\frac{2}{3}q_{t_j}^*}=  
\frac{1}{2 (q_{t_j}^{\ast})^2} .
\end{align*}
Because $\abs{x - \frac{p}{q}} < \frac{1}{q^2}$ implies $\frac{p}{q} = 
\frac{p_n}{q_n}$ for some $n\ge0$ and $0 < x - \frac{p}{q} < \frac{1}{q^2}$ 
implies $\frac{p}{q} = \frac{p_{m}^{\ast}}{q_{m}^{\ast}}$ for some 
$m \ge 0$, we have that all the $0$-convergents $\frac{p_{t_j}}{q_{t_j}}$ with $b_{{t_j}+1} \geq 4$ are even standard convergents, with possibly some $0$-convergents $\frac{p_{t_j}}{q_{t_j}}$ with $b_{{t_j}+1} =3$ in between, while the $0$-convergents with $b_{{t_j}+1} =2$ are not standard convergents. In conclusion, there exists a sequence $\{T_j\}\subset I_{**}(x)$ such that 
\[
\left\{
\begin{array}{ccc}
q_{T_j}^{\ast} & = & q_{2j} \\
a_{2j+1} & = & b_{T_j + 1} -2 
\end{array}
\right.
\]
for $j \ge 1$ (see Proposition \ref{prop1}). Clearly adding or omitting some indices $t_j$ with $b_{t_j+1}=3$ does not change the finiteness of the series:
\[
\abs{\sum_{j=0}^{\infty} 
\frac{\log(b_{t_j + 1})}{q_{t_j}^{\ast}} -
\sum_{j=1}^{\infty} \frac{\log a_{2j + 1}}{q_{2j}}} < 2\log 2 .
\]
Finally we have 
\[
\abs{\sum_{n=0}^{\infty} \frac{\log(b_{n+1}- 1)}{q_{n}^{\ast}} -
\sum_{j=1}^{\infty} \frac{\log a_{2j + 1}}{q_{2j}}} < 4\log 2 . 
\]
The same method holds for $0 < x < \frac{1}{2}$.
If $-x$ is a semi-Brjuno number, then so is $1 - x$.  By Lemma \ref{lemma1}, it is 
equivalent to 
\[
\sum_{j=0}^{\infty} \frac{\log a_{2j + 2}}{q_{2j+1}} < \infty. 
\]
Thus we have the assertion of the theorem. 
\end{proof}

\begin{rem}
We prefer not to extend the notion of the $(\alpha, u)$-Brjuno 
functions to the case $\alpha = 0$, because for general $u$, the divergence 
of $B_{(0,u)}(x)$ does not imply that $x$ has good by-excess  
approximations with respect to $u$ as in the case when $u(x) = \log x$. 
We have 
\begin{displaymath} 
\beta_{n}^{\ast} =  q_{n}^{\ast}x - p_{n}^{\ast} = \frac{x_n}{q_{n}^{\ast}  - q_{n-1}^{\ast} x_n} > \frac{1}{2} \frac{1}{(q_{n}^{\ast}  - q_{n-1}^{\ast} x_n)}  > \frac{1}{2q_{n}^{\ast}}
\end{displaymath} 
whenever $x_{n} > \frac{1}{2}$ and $b_{n} = 2$. 
Suppose that $u(x) = \frac{1}{\sqrt{x}}$, for example. 
Then it turns out that $\sum_{n \in I_{\ast}(x)} \frac{1}{q_{n}^{\ast}}=\infty$ 
implies $\sum_{n : x_n > 1/2} \beta_{n-1}^{\ast} u(x_n) = \infty$, 
which means that there exist good approximations from the \emph{right-hand side}, and not from the left-hand side. 
Moreover, the divergence is independent of this particular choice of $u$.  

\end{rem}

\end{document}